\documentclass[letterpaper, 10 pt, conference]{styles/ieeeconf}
\IEEEoverridecommandlockouts

\usepackage[utf8]{inputenc}
\usepackage[english]{babel}

\usepackage{cite}
\usepackage{amsmath,amssymb,amsfonts}
\usepackage{mathabx} 
\usepackage{xpatch}
\usepackage{dsfont}
\usepackage{acronym}
\usepackage{mathtools}
\usepackage{units}
\usepackage{comment}
\usepackage{thmtools}
\usepackage{graphicx}
\usepackage{algorithm}
\usepackage{hyperref}
\usepackage{algpseudocode}
\usepackage{url}            
\usepackage{booktabs}       
\usepackage{nicefrac}       
\usepackage{microtype}      
\usepackage{textcomp}
\usepackage[short]{optidef}
\usepackage{xcolor}
\definecolor{oxfordblue}{RGB}{0,33,71}          
\usepackage{subcaption}
\usepackage{caption}
\usepackage{tikz}
\usetikzlibrary{calc}
\usetikzlibrary{plotmarks}
\usetikzlibrary{positioning, fit}
\usepackage{pgfplots}
\usepgfplotslibrary{statistics, fillbetween, groupplots, external, colorbrewer, colormaps}
\pgfplotsset{compat=1.18}
\usepackage{xparse}

\newcommand\Prob[2]{\mathbb{P}_{#1}\Big\{ #2 \Big\}}

\DeclareRobustCommand\J[3]{
  {J_{#1}^{#2}(#3)} 
}
\DeclareRobustCommand\V[5]{\ifthenelse{\isempty{#4}}%
  {\ifthenelse{\isempty{#5}}
  {V^{\mu}_{#2, #3}}
  {V^{\mu}_{#2, #3}(#5)}}
  {\ifthenelse{\isempty{#5}}
  {V^{\mu}_{#2, #3, #4}}
  {V^{\mu}_{#2, #3,#4}( #5)}}
  }
\DeclareRobustCommand\Vb[5]{\ifthenelse{\isempty{#4}}%
  {\ifthenelse{\isempty{#5}}
  {V_{#2, #3}}
  {V_{#2, #3}(#5)}}
  {\ifthenelse{\isempty{#5}}
  {V_{#2, #3, #4}}
  {V_{#2, #3,#4}( #5)}}
  }
\DeclareRobustCommand\Measure[3]{\ifthenelse{\isempty{#2}}%
  {\ifthenelse{\isempty{#3}}
  {#1}
  {#1(#3)}}
  {\ifthenelse{\isempty{#3}}
  {#1_{#2}}
  {#1_{#2}({#3})}}
}
\renewcommand\P[2]{\Measure{\mathbb{P}}{#1}{#2}}

\newcommand\E[2]{
  {\mathbb{E}_{#1}\big[{#2}\big]}}

\NewDocumentCommand{\Risk}{ m m o o}{%
  \mathcal{R}_{#1}%
  \IfNoValueTF{#3}{}{^{#3}}%
  \big[{#2}%
  \IfNoValueTF{#4}{\big]}{}%
}
\newcommand\Unc[3]{
 \ifthenelse{\isempty{#3}}%
  {\mathcal{U}_{#1}\big[{#2}\big]}
  {\mathcal{U}_{#1}^{#3}\big[{#2}\big]}
  }

\newcommand\Opt[4]{
  \ifthenelse{\isempty{#2}}%
  {\mathrm{{#1}}}
  {\underset{\substack{\displaystyle #2}}{\mathrm{{#1}}}}
  \;
  \ifthenelse{\isempty{#4}}%
  {{#3}}
  {\Big\{ {#3}\  \Big| \  {#4} \Big\}}
}

\newtheorem{proposition}{Proposition}[section]

\newtheorem{lemma}{Lemma}[section]

\newtheorem{remark}{Remark}[section]
\newtheorem{assumption}{Assumption}[section]

\newcommand\numberthis{\addtocounter{equation}{1}\tag{\theequation}}
\title{On the suboptimality of stochastic MPC with varying constraint horizon }

\author{Allan Andre Do Nascimento, Andre Bertolace, Antonis Papachristodoulou and Kostas Margellos
\thanks{AAdN, AP and KM acknowledge funding support by MathWorks. AP was supported in part by UK’s Engineering and Physical Sciences Research Council projects EP/X017982/1, EP/Y014073/1 and UKRI2108. 
For the purpose of Open Access, the authors have applied a CC BY public copyright licence to any Author Accepted Manuscript (AAM) version arising from this submission. }
\thanks{All authors are with the Department of Engineering Science, University of Oxford,  Parks Road, Oxford OX1 3PJ, UK (email: \{allan.adn, andre.bertolace, antonis, kostas.margellos\}@eng.ox.ac.uk).} }

\begin{document}
\maketitle
\begin{abstract}
Enforcing stochastic state constraints over the full prediction horizon in Model Predictive Control (MPC) can be computationally demanding. Here we study stochastic MPC without terminal ingredients in which chance constraints are enforced only over a shorter constraint horizon. Using stochastic relaxed dynamic programming, we derive an explicit upper bound on the average expected closed-loop cost that depends on both prediction and constraint horizons. For linear quadratic problems with affine chance constraints and bounded uniform disturbances, we provide a deterministic reformulation via coordinate transformation and constraint tightening. Simulations illustrate the trade-off between computational effort and performance.
\end{abstract}
\vspace{-0.2cm}
\section{Introduction}

Modern control systems operate under state constraints and external disturbances. Model Predictive Control (MPC) enforces constraints by repeatedly solving a finite horizon constrained optimization problem, providing an approximation of the generally intractable optimal control problem (OCP) \cite{kouvaritakis2016model}. A central question is how close this approximation is to optimality, typically assessed through the closed loop or long run average cost. Performance guarantees depend on whether terminal ingredients are included. With terminal elements, the open loop cost is often used as an upper bound. Without them, guarantees are commonly expressed via suboptimality factors relating the finite horizon value function to the infinite horizon cost. The existence and form of these bounds depend on whether the setting is deterministic, stochastic, or distributionally robust, and on whether constraints are enforced over the full prediction horizon or only over a shorter constraint horizon.

For deterministic MPC with terminal ingredients, closed loop cost and stability guarantees are well established \cite{rawlings2020model,kouvaritakis2016model,grne2013nonlinear}. Similar results hold in the stochastic case when suitable terminal conditions are imposed \cite{rawlings2020model,kouvaritakis2016model}. Without terminal ingredients, deterministic MPC still admits performance guarantees under sufficiently long horizons or via relaxed dynamic programming (RDP) arguments \cite{grune2008infinite,grne2013nonlinear}. Deterministic schemes with variable constraint horizons have also been analyzed \cite{do2025constraint}. In contrast, stochastic and distributionally robust MPC without terminal ingredients and with variable constraint horizons lack comparable closed loop cost bounds, with existing results limited to fixed constraint horizons \cite{lorenzen2019stochastic}.

In this work, we address this gap by studying stochastic MPC with a variable constraint horizon and no terminal ingredients. Building on RDP ideas \cite{grne2013nonlinear} and related developments \cite{nascimento2025model}, we derive an explicit suboptimality bound that accounts for the distinct horizons employed. The results are developed for a parameterized stochastic MPC controller, which in the linear quadratic case can be solved deterministically via coordinate transformation and constraint tightening.

\vspace{-0.2cm}
\section{Problem formulation}
\label{sec:pf}

Consider the discrete-time system
\begin{equation}
\label{eq:sys}
    x_{k+1} = f(x_k, u_k, w_k),
\end{equation}
where $x_k \in \mathbb{R}^n$ is the state, $u_k \in \mathbb{R}^m$ the control input, and $w_k \in \mathcal{W}\subset \mathbb{R}^d$ encodes a disturbance vector. We assume that $\mathcal{W}$ is compact and consider admissible (non-anticipative) control policies.  
\begin{assumption}[Independent disturbances]
\label{assump:independence} The \textcolor{black}{random vectors} $W_k$, $k\geq0$, are independent and identically distributed (i.i.d.) with probability distribution $\P{}{}$ \textcolor{black}{on $\mathcal{W}\subseteq\mathbb{R}^d$}.
\end{assumption}
Under Assumption~\ref{assump:independence}, \eqref{eq:sys} satisfies the controlled Markov property. Hence, restricting attention to Markov policies entails no loss of optimality in the standard stochastic dynamic programming framework \cite{Bertsekas1995,Bertsekas1996}. For the subsequent derivations, we impose the following assumption.
\begin{assumption}
\label{ass1}
In System \eqref{eq:sys}, $f: \mathbb{R}^n \times \mathbb{R}^m \times \mathbb{R}^d \rightarrow \mathbb{R}^n$ is continuous with $f(0,0,0)=0$. 
\end{assumption}
We model predicted future states $x_{i+1|k}$ as random variables,
\begin{equation}
    x_{i+1|k} = f(x_{i|k}, u_{i|k}, W_{k+i}), \quad  W_{k+i} \sim \P{}{},
\end{equation}
with $x_{0|k} \overset{a.s.}{=} x_k$, i.e., predictions start \emph{almost surely} from the measured state at time $k$, $x_{i|k}$ is the open-loop prediction and $W_{k+i}$ is the disturbance random \textcolor{black}{vector}, both $i$ steps ahead of $k$, while $w_k$ is the disturbance realization at time $k$. We consider hard constraints on the controller,
\(
    u_{i|k} \in \mathcal{U},
\)
and probabilistic constraints on the state for all $k \ge0 $, 
\begin{equation}
\label{eq:prob_set_ctr}
    \Prob{}{x_{i+1|k} \in \mathcal{X} \mid x_{i|k}} \geq 1 - \varepsilon, \varepsilon \in \textcolor{black}{(0,1)}.
\end{equation}
where $1-\varepsilon$ is the constraint confidence level. Enforcing hard constraints on the control inputs is motivated by physical actuator limits, while chance constraints on the states can be introduced as a design choice to trade performance for ``controlled'' set violations. The cost function to be minimized is the expected value of a finite-horizon cost, i.e.,
\begin{equation}
    \J{N}{}{x_k, \mathbf{u}_{N|k}} = \E{\P{}{}}{\sum_{i=0}^{N-1} \ell(x_{i|k}, u_{i|k})\big|x_k},
\end{equation}
where \( \mathbf{u}_{N|k} = [u_{0|k}(x_k), \ldots, u_{N-1|k}(x_k)]^T\), are mappings of $x_k$, and the stage cost $\ell : \mathbb{R}^n \times \mathbb{R}^m \to \mathbb{R}_{\geq 0}$ is continuous and positive definite. The expected value operator $\mathbb{E}_{\P{}{}}$ is associated to $\P{}{}$.

Enforcing constraints at every prediction step can be expensive. We therefore study an MPC formulation in which the \emph{constraint horizon} differs from the prediction horizon. Although this setting can also be covered by \emph{multiple} and \emph{variable} constraint horizons \cite{nascimento2025model}, we focus on a partially constrained formulation, similar to \cite{do2025constraint}. \textcolor{black}{Here nonetheless, we focus on disturbed stochastic rather than nominal undisturbed systems, providing a lever to regulate constraint satisfaction while potentially reducing computational effort, in the disturbed setting. The stochastic suboptimality derivation extends the set-function-difference approach of \cite{nascimento2025model}, rather than the induction-based argument of \cite{do2025constraint}. Unlike \cite{nascimento2025model}, the second horizon portion is unconstrained and subject to no compactness requirement.} Relative to \cite{lorenzen2019stochastic}, we further allow partial constraint enforcement, \textcolor{black}{producing a suboptimality expression containing the prediction and constraint horizons, and derive a stochastic-robust ``quantile propagation" }for certain types of problems.

Specifically, probabilistic state constraints are enforced only over the constraint horizon $N-\tilde{N}$ of the prediction ($N-\tilde{N}+1$ predicted states in total), while the remaining $\tilde{N}-1$ terminal steps are left unconstrained. 

This yields a \emph{partially constrained} MPC problem, producing the value function $\Vb{\P{}{}}{N}{\tilde{N}}{}{x_k}$ and its minimizer, denoted by $\mathbf{u}_{N|k}^\star=[u^\star_{0|k}, \ldots, u^\star_{N-1|k}]^T$,
\begin{mini}
  {\mathbf{u}_{N|k}}
  {\J{N}{}{x_k, \mathbf{u}_{N|k}}} 
  {\label{prob:smpc}}
  {\Vb{\P{}{}}{N}{\tilde{N}}{}{x_k}:=}
  \addConstraint{x_{i+1|k}=f(x_{i|k}, u_{i|k}, W_{i+k})}
  \addConstraint{x_{0|k} = x_k}
  \addConstraint{u_{i|k} \in \mathcal{U},  \; i=0, \dots, N-1}
  \addConstraint{\Prob{}{x_{i+1|k} \in \mathcal{X} \mid x_{i|k}} \geq 1 - \varepsilon}
  \addConstraint{ i=0, \dots, N-\tilde{N}.}
\end{mini}
We consider the resulting sequence of control mappings, which are solutions of problem \eqref{prob:smpc} to be a parameterized (and generally suboptimal) class of control laws. In particular, we constrain the mappings $u^\star_{N-n|k}(\cdot)$ to admit a finite-dimensional parametrization, 
\begin{equation}\label{eq:paropt}
    u^\star_{N-n|k}(x_{N-n|k},v_{N-n|k}),\, n\in\{1,\dots,N\}
\end{equation}
where the closed loop depends on $x$, but $v_{N-n|k}\in\mathbb{R}^m$ \textcolor{black}{is a state-independent optimization parameter}. Note that under this parametrization, optimizing over \(u\) is equivalent to optimizing over the \(v\), since the control law is expressed in terms of \(v\). Hence, whenever we refer to optimizing over \(u\), it is implicitly understood that the optimization is carried out over \(v\). We adopt the following standard viability assumption used in MPC without terminal invariant set enforcement \cite{grne2013nonlinear,lorenzen2019stochastic}.
\begin{assumption}[Viability]
\label{ass2}
There exists a compact robust control invariant set \textcolor{black}{$\mathcal{X}_{\delta} \subseteq \mathcal{X}$}, such that for all $x_k \in \mathcal{X}_{\delta}$, $2 \leq \tilde{N} \leq N$ and $w \in \mathcal{W}$, $f(x_k,u^\star_{0|k},w) \in \mathcal{X}_{\delta}$ and $\Vb{\P{}{}}{N}{\tilde{N}}{}{}(x_k)$ is attained, for $\Vb{\P{}{}}{N}{\tilde{N}}{}{}(x_k)$ continuous on $\mathcal{X}_{\delta}$.
\end{assumption}
This assumption which affects only the \emph{first input} guarantees well-posedness and recursive feasibility of the MPC. The probabilistic constraint is enforced for the remainder of the constraint horizon $N-\tilde{N}$. The closed-loop controller is defined to be $u^\star_{0|k}(x_k)$. This means that at each instance $k$, \eqref{prob:smpc} is solved and $u^\star_{0|k}(x_k)$ is applied to the system
\begin{equation}\label{eq:closed-loop-dynamics}
    x_{k+1} = f(x_{k},u^\star_{0|k}(x_k),w_k).
\end{equation}
Problem \eqref{prob:smpc} is repeatedly solved along the closed-loop trajectory at states $x_0,\dots,x_k$, yielding closed loop inputs $u^\star_{0|0}(x_0),\dots,u^\star_{0|k}(x_k)$ that are applied to the plant. 
Under the realized disturbances $w_0,\dots,w_k$, this generates the closed-loop dynamics \eqref{eq:closed-loop-dynamics}. 
The associated \emph{average expected infinite-horizon closed-loop cost} is defined as
\begin{equation}\label{eq:closed-loop-cost}
   \bar{J}_{\infty}^{N,\tilde N}(x_0)= \limsup_{m\rightarrow \infty} \frac{1}{m}\sum^{m-1}_{k=0}\E{\P{}{}}{\ell(x^\star_{0|k}, u^\star_{0|k}(x_k))}.
\end{equation}
We study the average expected closed-loop optimal cost $\bar{J}_{\infty}^{N,\tilde N}(x)$ of \eqref{prob:smpc} under the control parameterization \eqref{eq:paropt}, and how the prediction horizon $N$ and constraint horizon $N-\tilde N$ affect a specific upper bound of \eqref{eq:closed-loop-cost}. This analysis yields performance bounds and can be used for stability analysis\cite{do2025constraint}. Note that we write $u^\star_{0|k}$ and $u^\star_{0|k}(x_k)$ interchangeably.

\vspace{-0.35cm}
\section{Main result}
\label{sec:fdc}
\subsection{Stochastic Closed Loop Sub-optimality}
\label{sclperf}

We begin by deriving a closed-loop upper bound for partially constrained stochastic MPC using stochastic RDP.

\begin{proposition}[Stochastic RDP] \label{prop:srdp}
    Consider $N\ge 2$, $2\le \tilde N\le N$. Assume there exists $c_2 \geq 0$ and $\alpha \in (0,1)$, such that, for all  $x_k \in \mathcal{X}_{\delta}$, with $x_0 \in \mathcal{X}_{\delta}$
        \begin{equation}\label{eq:srdp_1}
        \begin{aligned}
        & \E{\P{}{}}{\Vb{\P{}{}}{N}{\tilde N}{}{f(x_k,u^\star_{0|k}(x_k),w_k)}\big|x_k} - \Vb{\P{}{}}{N}{\tilde N}{}{x_k} \le \\
        & -\alpha \ell(x_k,u^\star_{0|k}(x_k)) + c_2.
        \end{aligned}    
        \end{equation} 
     Then, the average expected closed loop cost \eqref{eq:closed-loop-cost} incurred in running the closed loop dynamics \eqref{eq:closed-loop-dynamics} under MPC \eqref{prob:smpc} is
    \begin{equation} \label{ubcl}
    \bar{J}_{\infty}^{N,\tilde N}(x_0) \leq  \frac{c_2}{\alpha}. 
\end{equation}
\end{proposition}\vspace{2mm}
The proof can be found in Section \ref{appendix}. It follows similar steps as its deterministic counterpart \cite{do2025constraint}, with stochastic concepts adopted from \cite{lorenzen2019stochastic}. Now we propose explicit characterizations on $\alpha$ and $c_2$ based on the adopted control strategy. 

\subsection{Explicit expressions for $\alpha$ and $c_2$} \label{subs:Exp_alp_c}
Now, we derive different sequences of $\Vb{\P{}{}}{n}{\tilde N}{}{}{}(x_k),\ n\in{0,\ldots,N}$, defined backwards with a growing $n$. For $n\in\{1,\ldots,\tilde N -1\}$, the OCP's final part is
\begin{mini}
  {\mathbf{u}_{n|k}}
  {\E{\P{}{}}{\sum_{i=N-n}^{N-1} \ell(x_{i|k}, u_{i|k})\big|x_k}}
  {\label{eq:value-func}}
  {\Vb{\P{}{}}{n}{\tilde{N}}{}{x_k}:=}
  \addConstraint{x_{i+1|k} = f(x_{i|k}, u_{i|k}, W_{i+k})}
  \addConstraint{x_{N-n|k} = x_k}
  \addConstraint{u_{i|k} \in \mathcal{U}}
  \addConstraint{i=N-n, \dots, N-1.}
\end{mini}
where \( \mathbf{u}_{n|k} = [u_{N-n|k}, \ldots u_{N-1|k}]^T\), whereas, for $n\in\{\tilde N,\ldots,N\}$, the OCP's initial (chance constrained) part is
\begin{mini}
  {\mathbf{u}_{n|k}}
  {\E{\P{}{}}{\sum_{i=N-n}^{N-1} \ell(x_{i|k}, u_{i|k}) \big|x_k}}
  {\label{eq:value-func-2}}
  {\Vb{\P{}{}}{n}{\tilde{N}}{}{x_k}:=}
  \addConstraint{x_{i+1|k} = f(x_{i|k}, u_{i|k}, W_{i+k})}
  \addConstraint{x_{N-n|k} = x_k}
  \addConstraint{u_{i|k} \in \mathcal{U},\; i=N-n, \dots, N-1}
  \addConstraint{\P{}{}\{x_{i+1|k} \in \mathcal{X} \mid x_{i|k} \} \geq 1 - \varepsilon} 
  \addConstraint{i=N-n, \dots, N-\tilde{N}}.
\end{mini}
When $n=N$, \eqref{eq:value-func-2} is equivalent to \eqref{prob:smpc}. As no terminal cost is adopted in \eqref{prob:smpc}, for $n=0$ the value function $\Vb{\P{}{}}{0}{\tilde{N}}{}{x_k} := 0$. These OCP formulations arise naturally due to our variable-constraint horizon problem formulation. The probabilistic constraints \eqref{eq:prob_set_ctr} are incorporated in \eqref{eq:value-func-2} but not in \eqref{eq:value-func}. The use of a varying constraint horizon also reflects on our control formulation. For $n\in\{1,\ldots,\tilde N-1\}$, $\{u_n(x)\}_{n=1}^{\tilde N-1}$ and value functions $\{\Vb{\P{}{}}{n}{\tilde N}{}{x}\}_{n=1}^{\tilde N-1}$ satisfy
\begin{subequations}
\begin{equation}
        \Vb{\P{}{}}{n}{\tilde N}{}{x} :=  \min_{u \in \mathcal{U}} \left(\ell(x,u) +  \Vb{\mathbb{P}}{n-1}{\tilde N}{}{f(x,u,W)} \right) \label{eq:value-iteration3},
\end{equation}
and for $n\in\{\tilde N,\ldots,N\}$,
\begin{equation}
        \Vb{\P{}{}}{n}{\tilde N}{}{x} :=  \min_{u \in \mathcal{U}_n(x)}   \left(\ell(x,u) +  \Vb{\mathbb{P}}{n-1}{\tilde N}{}{f(x,u,W)} \right) \label{eq:value-iteration4},
\end{equation}
\end{subequations}
where $\mathcal{U}_n(x) = \{u \in \mathcal{U},\: \textrm{s.t.}, \: \mathbb{P}\{x_{N+1-n|k} \in \mathbb{X} \mid x_{N-n|k} \} \geq 1 - \varepsilon \}$. Unlike \eqref{eq:value-iteration3}, where the minimization is over $\mathcal{U}$, \eqref{eq:value-iteration4} minimizes over a subset further restricted by the input constraints induced by the chance constraints.
Explicit expressions for $\alpha$ and $c_2$ are obtained from the system, disturbance, and controller characteristics by partitioning the optimal open-loop state–input pair solving \eqref{prob:smpc} and decomposing the value function $\Vb{\P{}{}}{N}{\tilde{N}}{}{x_k}$ as,
\begin{subequations}
    \begin{equation}
 \label{eq:brkdwn}
    \Vb{\P{}{}}{N}{\tilde{N}}{}{x_k} = \sum^{N-\tilde{N}}_{n=0} \lambda_{s_1}[(n,k)|x_k] + \sum^{N-1}_{n=N-\tilde{N}+1} \lambda_{s_2}[(n,k)|x_k],
    \end{equation}
where
\begin{equation}
    \lambda_{s_1}[(n,k)|x_k] = \E{\P{}{}}{\ell(x^\star_{n|k},u^\star_{n|k})|x_k}{},
\end{equation}
for $n=0,\dots,N-\tilde{N}$, and 
\begin{equation}
    \lambda_{s_2}[(n,k)|x_k] =  \E{\P{}{}}{\ell(x^\star_{n|k},u^\star_{n|k})|x_k}{},
\end{equation}
for $n=N-\tilde{N}+1,\dots,N-1$.
\end{subequations}

In \eqref{eq:brkdwn}, \( \Vb{\P{}{}}{N}{\tilde{N}}{}{x_k} \) is decomposed to separate the running cost contributions from the constrained and unconstrained portions of \eqref{prob:smpc}. Here, \(\lambda_{s_1}[(n,k)|x_k]\) corresponds to predicted states and inputs satisfying the probabilistic constraints, while \(\lambda_{s_2}[(n,k)|x_k]\) represents the remaining \textcolor{black}{predicted states and inputs}, with $[(n,k)|x_k]$ denoting the running cost of the state and input predicted $n$ steps ahead of time $k$. \textcolor{black}{We highlight that as in \cite{lorenzen2019stochastic}, when $V_{n,\tilde{N}}(x_{j\mid k})$ is evaluated at the predicted random state $x_{j\mid k}$, the expectation in its definition remains conditioned on $x_k$}. Consider Lemma \ref{lem:eqVNNtilde}, important for the main result.


\begin{lemma}\label{lem:eqVNNtilde}
    For $N \geq 3, 2 \leq \tilde{N} \leq N-1$, consider
\begin{eqnarray*}
  \Vb{\P{}{}}{N}{\tilde{N}}{}{x_{k+1}} &=& \sum^{N-\tilde{N}}_{n=0} \lambda_{s_1}[(n,k+1)|x_{k+1}] \nonumber \\ &&+ \sum^{N-1}_{n=N-\tilde{N}+1} \lambda_{s_2}[(n,k+1)|x_{k+1}]. \notag
  \end{eqnarray*}
   Then, the following upper bound holds:
\begin{align} \label{eq:ubv}
        \textcolor{black}{\E{\P{}{}}{\Vb{\P{}{}}{N}{\tilde{N}}{}{x_{k+1}}\big| x_k}} & \leq \sum^{N-\tilde{N}}_{n=1} \lambda_{s_1} [(n,k)|x_{k}]  \nonumber \\
        & + \Vb{\P{}{}}{\tilde{N}}{\tilde{N}}{}{x^\star_{N-\tilde{N}+1|k}}.
\end{align}
\end{lemma}
The proof (in Section \ref{appendix}) follows arguments in \cite{nascimento2025model,lorenzen2019stochastic}. We now use Lemma~\ref{lem:eqVNNtilde} to characterize the left-hand side of \eqref{eq:srdp_1}:
\begin{equation}
\begin{aligned}
    \E{\P{}{}}{\Vb{\P{}{}}{N}{\tilde N}{}{x_{k+1}} \big| x_k}- \Vb{\P{}{}}{N}{\tilde N}{}{x_k} \leq {\sum^{N-\tilde{N}}_{n=1} \lambda_{s_1}[(n,k)|x_{k}]} \\
     + {\Vb{\P{}{}}{\tilde N}{\tilde N}{}{x^\star_{N-\tilde{N}+1|k}}} - \E{\P{}{}}{\ell(x_{0|k}, u_{0|k})} \\
     - {\sum^{N-\tilde{N}}_{n=1} \lambda_{s_1}[(n,k)|x_{k}]} -{\sum^{N-1}_{n=N-\tilde{N}+1}\lambda_{s_2}[(n,k)|x_k]} 
    \label{eq:exp_ub}
\end{aligned}    
\end{equation}
Discarding $- \sum^{N-1}_{n=N-\tilde{N}+1}\lambda_{s_2}[(n,k)|x_k] \leq 0$ and simplifying \eqref{eq:exp_ub}: 
\begin{equation}
\begin{aligned}
    & \E{\P{}{}}{\Vb{\P{}{}}{N}{\tilde N}{}{x_{k+1}} \big| x_k}- \Vb{\P{}{}}{N}{\tilde N}{}{x_k} \leq\\
    & - \E{\P{}{}}{\ell(x_{0|k}, u_{0|k})} + {\Vb{\P{}{}}{\tilde{N}}{\tilde N}{}{x^\star_{N-\tilde{N}+1|k}}}. 
    \label{eq:exp_ub_simp}
\end{aligned}
\end{equation}
Based on \eqref{eq:exp_ub_simp}, we obtain an explicit expression for $\alpha$, by adopting a modified version of cost controllability \cite{nascimento2025model}:

\begin{assumption}[Cost controllability] \label{ass:stctrl}

We assume that \textcolor{black}{for each fixed pair $(N,\tilde N)$, there exist constants $C_1>0$ and $\sigma_1\in(0,1)$, independent of $x_k$, such that, for all $x_k\in\mathcal{X}_{\delta}$,} the optimal running cost at step $N-\tilde{N}$ is \emph{cost controllable}:
\begin{equation}
\lambda_{s_1}[(N-\tilde{N},k)|x_k]
\le
C_1 \sigma_1^{N-\tilde N}\ell(x^\star_{0|k},u^\star_{0|k})
+
\ell^0_{N-\tilde N},
\label{eq:stctrl1}
\end{equation}
meaning the expected value of the running cost at the ``switching step'' is cost bounded. Moreover, there exist constants $C_2>0$ and $\sigma_2\in(0,1)$, \textcolor{black}{independent of $x_k$, such that, for all $x_k\in\mathcal{X}_{\delta}$, there exists} a state/control pair
$\{\tilde x_{N-\tilde N+1+j|k},\tilde u_{N-\tilde N+1+j|k}\}_{j=0}^{\tilde N-1}$
admissible for the corresponding $\tilde N$-stage tail problem, starting at
$\tilde x_{N-\tilde N+1|k}=x^\star_{N-\tilde N+1|k}$, such that:
\begin{align*}
& \mathbb E_{\mathbb P}\!\left[
\ell(\tilde x_{N-\tilde N+1+j|k},
\tilde u_{N-\tilde N+1+j|k})
\,\middle|\,x_k
\right] \leq \\
& C_2 \sigma_1 \sigma_2^{j}
\lambda_{s_1}[(N-\tilde{N},k)|x_k]
+\ell^0_{N-\tilde N+1+j}
+\Delta^0_{N-\tilde N}
\label{eq:stctrl2}\numberthis \\
& j=0,\dots,\tilde N-1, \\
& \textrm{where} \\
& \ell_i^0
=
\mathbb E_{\mathbb P}\!\left[
\ell(x^\star_{i|k},u^\star_{i|k})
\,\middle|\,x_k
\right],
\text{ when }x^\star_{0|k}=0,
\end{align*}
\textcolor{black}{$i=0,\ldots,N-\tilde N$,
\begin{align*}
\ell^0_{N-\tilde N+1+j}
=
\mathbb E_{\mathbb P}\!\left[
\ell(\tilde x^\star_{N-\tilde N+1+j|k},
\tilde u^\star_{N-\tilde N+1+j|k})
\,\middle|\,x_k
\right],
\end{align*}
$j=0,\ldots,\tilde N-1$, where
$(\tilde x^\star,\tilde u^\star)$ denotes an optimizer of the corresponding auxiliary problem
$V_{\tilde N,\tilde N}(x^\star_{N-\tilde N+1|k})$
when $x^\star_{0|k}=0$, and}
\begin{equation*}
\Delta^0_{N-\tilde N}
=
\frac{1}{\tilde N}
\sum_{r=0}^{N-\tilde N-1}\ell_r^0.
\end{equation*}
\end{assumption}

Bound \eqref{eq:stctrl1} captures the chance-constrained part decay with rate $\sigma_1$. Here we use a milder condition, linking $\lambda_{s_1}[(N-\tilde{N},k)|x_k]$ to $\ell(x^\star_{0|k},u^\star_{0|k})$, and its drift term instead of imposing it over the full constrained open-loop trajectory. Bound \eqref{eq:stctrl2} is a shifted version of stochastic cost controllability \cite{lorenzen2019stochastic}, accounting first for the chance-constrained term via $\sigma_1$ and then for the unconstrained states via $\sigma_2$.
\begin{remark}
    As noted in \cite{lorenzen2019stochastic}, any summable $\mathcal{KL}$ function could be used in \eqref{eq:stctrl2}. Nominal and stochastic cost controllability are common assumptions in nominal and stochastic MPC without terminal elements \cite{grne2013nonlinear, bold2024data, moldenhauer2025robust, lorenzen2019stochastic}. Here, we adopt a stronger version of cost controllability, to obtain a closed form sub-optimality estimate, potentially yielding tighter bounds \cite{lorenzen2019stochastic} (to be further explored in the future). \textcolor{black}{Given the requirements of Assumption \ref{ass:stctrl}, we expect its constants to be more simply computable or bounded a priori for exponentially stabilizable nonlinear systems admitting a suitable Lyapunov function, as well as for linear-quadratic systems with an admissible stabilizing feedback.}
\end{remark}

Now using Assumption \ref{ass:stctrl}, we revisit \eqref{eq:exp_ub_simp}, and characterize $\alpha$ explicitly via the following proposition:
\begin{proposition}
\label{propalpha} Consider Assumptions \ref{ass1}, \ref{ass2} and \ref{ass:stctrl}. Then, for $N \geq 3$, $2 \leq \tilde{N} \leq N-1$, and all $x_k\in\mathcal{X}_\delta$, Proposition \ref{prop:srdp} admits the following explicit expressions for $\alpha$ and $c_2$, \textcolor{black}{provided that the resulting $\alpha\in(0,1)$} -- Proof in Section \ref{appendix}: 
    \begin{equation}
    \label{prop:alpha_exp}
    \begin{aligned}
        & \alpha =  1- C_1C_2\sigma^{N-\tilde{N}+1}_1 \left ( \frac{1-\sigma^{\tilde{N}}_2}{1-\sigma_2} \right), \\
        & c_2 = \max \left \{1,C_2\sigma_1\frac{1-\sigma^{\tilde N}_2}{1-\sigma_2} \right \} \left (\sum^{N}_{i=0} \ell^0_i \right). 
    \end{aligned}
    \end{equation} 
\end{proposition}
This is an explicit upper bound on the average expected closed-loop cost \eqref{eq:closed-loop-cost}. Disturbances enter through $\ell_i^0$, the constrained portion through $\sigma_1^{N-\tilde N+1}$, and the tail through $\sum_{j=0}^{\tilde N-1}\sigma_2^j$. \textcolor{black}{With other quantities fixed, decreasing $\tilde N$ tightens the bound but adds chance constraints. Since these quantities depend on $(N,\tilde N)$, monotonicity is not guaranteed. The risk level $\epsilon$ enters through the feasible sets, hence through $\alpha$ and $c_2$.} As usual in RDP, the bound is meaningful only for $\alpha\in(0,1)$.
\begin{remark}
    \textcolor{black}{The bound obtained in \eqref{ubcl} can be conservative. This is due to the discarding performed in \eqref{eq:exp_ub}, uniform estimations of Assumption \ref{ass:stctrl}, and the $\max\{\cdot,\cdot\}$ operator in \eqref{prop:alpha_exp}. Tighter admissible constants/ decay rates and trajectory-wise (instead of set uniform) usage of Assumption \ref{ass:stctrl}  can reduce this conservatism.} Nonetheless, $\alpha$ can still be a useful design tool. When $\ell(\cdot,\cdot) > 0, (x,u) \neq 0$ and $\sigma$'s are predefined, \textcolor{black}{enforcing $\alpha>0$ via $(N, \tilde N)$ makes the expected drift in \eqref{eq:srdp_1} negative whenever $\ell(x_k,u^\star_{0|k}(x_k))>c_2/\alpha$. This provides the strict-drift condition required by stochastic stability arguments (\cite{lorenzen2019stochastic}, Proposition~1 and Corollary~1).}
\end{remark}

\subsection{Reformulation of Partially Constrained MPC}
Restricting \eqref{prob:smpc} to a class of linear systems enables its deterministic reformulation \cite{lorenzen2019stochastic}. We adopt: 
\begin{mini}
  {\mathbf{u}_{N|k}}
  {\mathbb{E}[\sum^{N-1}_{i=0} x^{\top}_{i|k}Qx_{i|k} + u^{\top}_{i|k}Ru_{i|k}\mid x_k]} 
  {\label{prob:smpclin}}
  {}
  \addConstraint{x_{i+1|k}=A x_{i|k}+B u_{i|k}+G W_{k+i}}
  \addConstraint{x_{0|k} = x_k}
  \addConstraint{u_{i|k} \in \mathcal{U},  \; i=0, \dots, N-1}
  \addConstraint{\Prob{}{g(x_{i+1|k},x_{i|k}) \le 0 \mid x_{i|k}} \geq 1 - \varepsilon}
  \addConstraint{ i=0, \dots, N-\tilde{N}.}
\end{mini}
Where $Q,R \succ 0$ and $g:\mathbb{R}^{n} \times \mathbb{R}^{n} \to \mathbb{R}$ is assumed to be affine, representing individual chance constraints. \textcolor{black}{The preceding analysis applies directly to conditional
chance constraints of the form imposed in \eqref{prob:smpclin}. Furthermore, for every state $x_{i|k}\in\mathcal{X}$ at which the
chance constraint is imposed, we assume that
$\{x_{i+1|k}\in\mathbb{R}^{n}:
g(x_{i+1|k},x_{i|k})\leq0\}\subseteq\mathcal{X}$,
so that the chance constraint in \eqref{prob:smpclin} is a sufficient
condition for \eqref{eq:prob_set_ctr}
whenever $x_{i|k}\in\mathcal{X}$.} The disturbance $W_{k+i}$ is assumed uniformly distributed over the bounded box $\mathbf{W}:=\{-w_{\max}\le W_{k+i}\le w_{\max}\}$, with independent components and zero mean. We adopt the parameterization \eqref{eq:paropt}, $u_{i|k}=v_{i|k}+Ke_{i|k}$, where $K\in\mathbb{R}^{m\times n}$ is chosen such that $A_{\mathrm{cl}}\coloneqq A+BK$ is Schur, and $\{v_{i|k}\}_{i=0}^{N-1}$ are the decision variables. We define the conditional mean and deviation:
\begin{equation}\label{eq:ze_def_app}
z_{i|k}\coloneqq \mathbb{E}[x_{i|k}\mid x_k], \,\,
e_{i|k}\coloneqq x_{i|k}-z_{i|k}, \,\,
x_{i|k}=z_{i|k}+e_{i|k},
\end{equation}
so that $e_{0|k}=0$. Using \eqref{eq:ze_def_app}, the original linear system can be decomposed into $z_{i+1|k} = A z_{i|k}+B v_{i|k},
\: z_{0|k}=x_k \textrm{ and } e_{i+1|k} = A_{\mathrm{cl}} e_{i|k}+G W_{k+i}, \: e_{0|k}=0$, nominal and deviation dynamics respectively. This allows us to reformulate \eqref{prob:smpclin} as a deterministic MPC problem by transforming the cost and tightening the constraints.

\subsubsection{Cost function transformation} 
The quadratic finite-horizon cost in \eqref{prob:smpclin} is split into a nominal cost and a disturbance propagation as
\begin{equation}
\begin{aligned}
\J{N}{}{x_k, \mathbf{u}_{N|k}} ={}&
\sum_{i=0}^{N-1}\left(z_{i|k}^\top Q z_{i|k}
+v_{i|k}^\top R v_{i|k}\right)
\\
&+\sum_{i=0}^{N-1}\left(\mathrm{tr}(Q \Sigma_{e,i})
+\mathrm{tr}(K^\top R K\,\Sigma_{e,i})\right),
\end{aligned}
\end{equation}
where $\Sigma_{e,i}\in\mathbb{R}^{n\times n}$ is the covariance matrix for the deviation, being propagated via the recursion at prediction step $i$ as $\textcolor{black}{\Sigma_W:=\mathrm{cov}(W_k)}, \: \Sigma_{e,0}:=0_{n \times n}, \Sigma_{e,i+1}:=A_{\mathrm{cl}}\Sigma_{e,i}A_{\mathrm{cl}}^\top+G\textcolor{black}{\Sigma_W} G^\top, \textrm{ for } i=0,\dots,N-1$.

\subsubsection{Input Constraint Tightening}
To guarantee $u_{i|k} \in \mathcal{U}$ for all admissible disturbances, we construct a robust tube around the nominal trajectory by deriving a sequence of deviation reachable sets. Specifically, we compute $\{E_i\}_{i=0}^{N}$, $E_i \subset \mathbb{R}^n$, using the Minkowski recursion,
\(
E_0:=\{0\},\: E_{i+1}:=A_{\mathrm{cl}}E_i\oplus G \mathbf{W}\,i=0,\dots,N-1,
\)
where $\oplus$ denotes the Minkowski sum of two sets defined by \cite{li2014sweep},
\(
M \oplus N \coloneq \{m + n \mid m \in M, n \in N\}.
\)
Since $e_{i|k}\in E_i$ for all admissible disturbance realizations, we enforce $
u_{i|k} = v_{i|k}+Ke_{i|k}\in \mathcal{U}, \: \forall e_{i|k}\in E_i$ (which is the assumed form of \eqref{eq:paropt}, given the coordinate transformation in \eqref{eq:ze_def_app}), by restricting the nominal inputs $v_{i|k}$ to the tightened sets
\(
V_i:=U\ominus (K E_i)
=\{v\in\mathbb{R}^{m}: v+(K E_i)\subseteq U\},\,i=0,\dots,N-1,
\)
where $\ominus$ denotes the Pontryagin set difference, this implies $v_{i|k}+Ke_{i|k}\in \mathcal{U}$. 

\subsubsection{State Constraint Tightening}
These constraints require propagation of \textcolor{black}{the disturbance set through the system dynamics, computation of the} disturbance distributions affecting the system and inversion of its Cumulative Distribution Functions (CDF). \textcolor{black}{Since \(g\) is affine, let \(a,b\in\mathbb{R}^n\) and \(c\in\mathbb{R}\) denote its coefficients. Then $g$ can be expressed as $g(x_{i+1|k},x_{i|k}) =a^\top x_{i+1|k}+b^\top x_{i|k}+c$. Using the error dynamics, we decompose the affine constraint into a deterministic term depending on $(z_{i+1|k},z_{i|k})$ and a scalar residual depending on the error},
\(
g(x_{i+1|k},x_{i|k}) = g(z_{i+1|k},z_{i|k}) + \delta_i.
\)
\textcolor{black}{$\delta_i =(a^\top A_{\mathrm{cl}}+b^\top)e_{i|k}
+a^\top G W_{k+i}$}.
\textcolor{black}{Here, \(c\) remains in \(g(z_{i+1|k},z_{i|k})\) and therefore does not appear in \(\delta_i\). For any scalar random variable \(Y\), we use $F^{-1}_{Y}(p), p\in(0,1)$, meaning the subscript \(Y\) identifies the random variable whose CDF is inverted. Note that, for a given \(x_k\), \(z_{i|k}\) is fixed, thus, conditioning on \(x_{i|k}\) fixes \(e_{i|k}\), and only \(W_{k+i}\) remains random. The chance constraint in \eqref{prob:smpclin} is therefore enforced for every \(e_{i|k}\in E_i\) by} $g(z_{i+1|k},z_{i|k}) + q_i\le 0$, where
\[
\textcolor{black}{
q_i:=
\max_{e\in E_i}(a^\top A_{\mathrm{cl}}+b^\top)e
+F^{-1}_{a^\top G W_{k+i}}(1-\varepsilon),
}
\]
\textcolor{black}{ for $i=0,\dots,N-\tilde N$. The first term accounts for the disturbance set propagated through the closed loop error dynamics, whereas the second is the \((1-\varepsilon)\)-quantile of the new-disturbance term. As the error dynamics and \(E_i\) are independent of the decision variables, so is \(q_i\).} Since we assume that \textcolor{black}{$a^\top G W_{k+i}$} is a weighted sum of independent uniform random variables, $q_i$ can be computed exactly. \textcolor{black}{The CDF of \(a^\top G W_{k+i}\)} is evaluated via the generalized Irwin--Hall distribution for independently  distributed uniform random variables \cite{sadooghi2009distribution, marengo2017geometric}. Then its inverse CDF \textcolor{black}{$F^{-1}_{a^\top G W_{k+i}}(\cdot)$} is evaluated at its confidence threshold $(1-\varepsilon)$.

\textcolor{black}{Alternatively, if the chance constraint is conditioned on \(x_k\) rather than \(x_{i|k}\), \(e_{i|k}\) also remains random and the complete-residual quantile gives the exact tightening}
\[
\textcolor{black}{
g(z_{i+1|k},z_{i|k})
+F^{-1}_{\delta_i}(1-\varepsilon)
\leq0.
}
\]
\textcolor{black}{However, the shifted candidate may then be infeasible, requiring the additional analysis in \cite{lorenzen2019stochastic}. Assumption \ref{assump:independence} imposes independence across time, but allows dependence within $W_k$, and any $\mathbb{P}$ satisfying it can be used in the main results. Uniformity and componentwise independence are used only for the exact quantile evaluation. Ambiguity sets would require time-consistent worst-case counterparts of some expressions. Distribution mismatch and temporal correlation are covered in upcoming works.}


\subsubsection{Corresponding deterministic MPC} \label{sbsc:cdmpc}

Equation \eqref{prob:smpclin} is then transformed into: 
\begin{mini}
  {\mathbf{v}_{N|k}}
  {%
    \sum_{i=0}^{N-1}\big(z_{i|k}^\top Q z_{i|k}
    +v_{i|k}^\top R v_{i|k}\big)
  }{\label{prob:tightmpc}}{}
  \breakObjective{%
    \;+\;\sum_{i=0}^{N-1}\Big(\mathrm{tr}(Q\,\Sigma_{e,i})
    +\mathrm{tr}(K^\top R K\,\Sigma_{e,i})\Big)
  }
  \addConstraint{z_{i+1|k}=A z_{i|k}+B v_{i|k}}
  \addConstraint{z_{0|k} = x_k}
  \addConstraint{v_{i|k} \in V_i,  \; i=0, \dots, N-1}
  \addConstraint{g(z_{i+1|k},z_{i|k}) + q_i\le 0}
  \addConstraint{i=0, \dots, N-\tilde{N}.}
\end{mini}

In this setting, the constant $c_2$ in Proposition~\ref{propalpha} becomes:
\begin{subequations}
\begin{align}\label{eq:c2}
&\max\left\{
1,C_2\sigma_1\frac{1-\sigma^{\tilde N}_2}{1-\sigma_2}
\right\}
\sum_{i=0}^{N}\Big(
\mathrm{tr}(Q\,\Sigma_{e,i})
+\mathrm{tr}(K^\top R K\,\Sigma_{e,i})
\Big),
\end{align}
while $\alpha$ is obtained by evaluating Assumption \ref{ass:stctrl} on the nominal part of \eqref{prob:tightmpc}:
\begin{mini}
  {\mathbf{v}_{N|k}}
  {%
    \sum_{i=0}^{N-1}\big(z_{i|k}^\top Q z_{i|k}
    +v_{i|k}^\top R v_{i|k}\big)
  }{\label{prob:nomtightmpc}}{}
  \addConstraint{z_{i+1|k}=A z_{i|k}+B v_{i|k}}
  \addConstraint{z_{0|k} = x_k}
  \addConstraint{v_{i|k} \in V_i,  \; i=0, \dots, N-1}
  \addConstraint{g(z_{i+1|k},z_{i|k}) + q_i\le 0}
  \addConstraint{i=0, \dots, N-\tilde{N}.}
\end{mini}
Together with \textcolor{black}{$V_{\tilde{N},\tilde{N}}(z^\star_{N-\tilde{N}+1|k})$, this problem yields the parameter $C_1$, $C_2$, $\sigma_1$, and $\sigma_2$. We will illustrate next Proposition \ref{propalpha} used with the reformulation just presented}.
\textcolor{black}{We assume that $0\in V_i$ at every relevant input stage and $g(0,0)+q_i\leq0$ for every tightened chance constraint. Since $Q,R\succ0$, the corresponding nominal optimizer initialized at zero is then identically zero.}
\end{subequations}

\input{sections/simulations.tex}
\section{Conclusion}
\label{sec:c}
This work analyzed the effect of partial constraint enforcement in stochastic MPC without terminal ingredients on the closed-loop upper bound. A parameterized finite-dimensional SMPC formulation was derived, and chance constraints were addressed via deterministic tightening, leading to a deterministic MPC reformulation. \textcolor{black}{Simulations illustrated the trajectory-based a posteriori $\widehat{\alpha}$-indicator and evaluated computational effort.}


\bibliographystyle{IEEEtran}
\bibliography{refs.bib}
\section{Appendix} \label{appendix}
\vspace{-0.2cm}
\subsection{Proof of Proposition \ref{prop:srdp}}
Since \eqref{eq:srdp_1} holds for all $x_k \in \mathcal{X}_{\delta}$, it also holds for $x^\star_{0|k}$. Rearranging terms gives, 
\begin{equation} 
\label{ublalpgen}
\begin{aligned}
    & \alpha \E{\P{}{}}{\ell(x^\star_{0|k}, u^\star_{0|k})} \\
    & \leq \Vb{\P{}{}}{N}{\tilde N}{}{x_k} - \E{\P{}{}}{\Vb{\P{}{}}{N}{\tilde N}{}{x_{k+1}} \big|x_k} + c_2
\end{aligned}
\end{equation}
Summing \eqref{ublalpgen} over m sequential time steps gives us:
\begin{equation}
\label{ubseqgen}
\begin{aligned}
    & \alpha \sum^{m-1}_{k=0}\E{\P{}{}}{\ell(x^\star_{0|k}, u^\star_{0|k})} \\
    & \leq \sum^{m-1}_{k=0} \left( \Vb{\P{}{}}{N}{\tilde N}{}{x_k} - \E{\mathbb{P}}{\Vb{\P{}{}}{N}{\tilde N}{}{x_{k+1}} \big|x_k}\right) + mc_2. 
\end{aligned}
\end{equation} 
Taking the conditional expectation of \eqref{ubseqgen} given $x_0$, i.e.,
\(
\E{\P{}{}}{\E{\P{}{}}{\cdot \mid x_k} \mid x_0},
\)
and applying the tower property of conditional expectation,
\(
\E{\P{}{}}{\E{\P{}{}}{Z \mid x_k} \mid x_0} = \E{\P{}{}}{ Z \mid x_0},
\)
we obtain an inequality conditioned on the initial state. Telescoping the sum on the right-hand side and dividing both sides by $m$ and $\alpha$ yields:
\begin{equation}
\label{ubseq}
\begin{aligned}
    & \frac{1}{m}\sum^{m-1}_{k=0}\E{\P{}{}}{\ell(x^\star_{0|k}, u^\star_{0|k})} \\
    & \leq \frac{ \E{\P{}{}}{\Vb{\P{}{}}{N}{\tilde N}{}{x_0} \big| x_0} - \E{\P{}{}}{\Vb{\P{}{}}{N}{\tilde N}{}{x_{m}}\big|x_0}}{m\alpha} + \frac{c_2}{\alpha}. 
\end{aligned}
\end{equation}
From Assumption \ref{ass2}, continuity of $\Vb{\P{}{}}{N}{\tilde N}{}{x_n}$ in the compact set $\mathcal{X}_{\delta}$ implies $\sup_n \E{\P{}{}}{\Vb{\P{}{}}{N}{\tilde N}{}{x_n}|x_0} < \infty$. Letting $m \to \infty$, yields:
\begin{equation}
\label{ubclgen}
    \bar{J}_{\infty}^{N,\tilde N}(x_0) = \limsup_{m\rightarrow \infty} \frac{1}{m}\sum^{m-1}_{k=0}\E{\P{}{}}{\ell(x^\star_{0|k}, u^\star_{0|k})} \leq  \frac{c_2}{\alpha}. 
\end{equation}
\vspace{-0.4cm}
\subsection{Proof of Lemma \ref{lem:eqVNNtilde}}
Using the breakdown proposed in \eqref{eq:brkdwn}:
\allowdisplaybreaks
\begin{equation}
\label{eq:exp_diff}
\begin{aligned}
    & \E{\P{}{}}{\Vb{\P{}{}}{N}{\tilde N}{}{x_{k+1}}\big|x_k} \\
    & = \E{\P{}{}}{\sum^{N-\tilde{N}-1}_{n=0} \lambda_{s_1}[(n,k+1)|x_{k+1}] \big| x_k} \\
    & + \E{\P{}{}}{\lambda_{s_1}[(N-\tilde{N},k+1)|x_{k+1}] \big| x_k} \\
    & + \E{\P{}{}}{\sum^{N-1}_{n=N-\tilde{N}+1} \lambda_{s_2}[(n,k+1)|x_{k+1}] \big| x_k}. 
\end{aligned}
\end{equation}
After the equality, the second term is a running cost with state and input under the chance constraint, and the third term is a sum corresponding to inputs not affected by it. Using \eqref{eq:exp_diff}, we upper bound $\E{\P{}{}}{\Vb{\P{}{}}{N}{\tilde N}{}{x_{k+1}}\big|x_k}$ by:
\begin{equation}\label{eq:exp_ubap}
    \begin{aligned}
    \E{\P{}{}}{\Vb{\P{}{}}{N}{\tilde N}{}{x_{k+1}}\big|x_k} & \leq \E{\P{}{}}{\sum^{N-\tilde{N}}_{n=1} \lambda_{s_1}[(n,k)|x_{k+1}] \big| x_k} \\ 
    & + \E{\P{}{}}{\Vb{\P{}{}}{\tilde N}{\tilde N}{}{x^\star_{N-\tilde{N}+1|k}}\big|x_k} \\
    = \sum^{N-\tilde{N}}_{n=1} \lambda_{s_1}& [(n,k)|x_{k}] + \Vb{\P{}{}}{\tilde{N}}{\tilde{N}}{}{x^\star_{N-\tilde{N}+1|k}}. 
   \end{aligned}
\end{equation}
The first summation in \eqref{eq:exp_diff} is replaced by the upper bound using $x^\star_{0|k+1} \overset{a.s.}{=} x^\star_{1|k}$ and the shifted suboptimal (with respect to $k+1$) inputs $[u^\star_{1|k}(x), \ldots u^\star_{N-\tilde{N}|k}(x)]^T$. The last two terms are upper-bounded by a value function of length $\tilde N$, with its first input constrained, starting at $x_{N-\tilde{N}+1|k}$. The inner expectation simplifies by the tower property. The bound in \eqref{eq:exp_ubap} follows by comparing the value function on the left with selected suboptimal terms on the right.
\subsection{Proof of Proposition \ref{propalpha}}
Using Assumption \ref{ass:stctrl}, on the right-hand side of \eqref{eq:exp_ub_simp}:
\allowdisplaybreaks
\begin{align*} \label{ub1}
    & \E{\P{}{}}{\Vb{\P{}{}}{N}{\tilde N}{}{x_{k+1}} \big |x_k}-\Vb{\P{}{}}{N}{\tilde N}{}{x_k} \\
    & \leq -\E{\P{}{}}{\ell(x_{0|k}, u_{0|k})} + {\Vb{\P{}{}}{\tilde N}{\tilde N}{}{x^\star_{N-\tilde{N}+1|k}}}\\
    &\leq  -\E{\P{}{}}{\ell(x_{0|k}, u_{0|k})} + \tilde{N}\Delta^{0}_{N-\tilde N} + \sum^{\tilde N}_{j=1}\ell^0_{N-\tilde{N}+j} \\
    & + C_2 \sigma_1 \left( \sum^{\tilde N -1}_{j=0} \sigma^j_2 \: \E{\P{}{}}{\ell(x^\star_{N-\tilde{N}|k}, u^\star_{N-\tilde{N}|k})}\right)  \\
    &\leq -\E{\P{}{}}{\ell(x_{0|k}, u_{0|k})} + \tilde{N}\Delta^{0}_{N-\tilde N} + \sum^{\tilde N}_{j=1}\ell^0_{N-\tilde{N}+j}\\
    & + C_1 C_2 \sigma^{N-\tilde{N}+1}_1 \left(\frac{1-\sigma^{\tilde N}_2}{1-\sigma_2} \right) \E{\P{}{}}{\ell(x^\star_{0|k}, u^\star_{0|k})}  \\
    & + C_2 \sigma_1 \frac{1-\sigma^{\tilde N}_2}{1-\sigma_2} \ell^0_{N-\tilde{N}} \numberthis. 
\end{align*}
For the last inequality, we used \eqref{eq:stctrl1}. Rearranging the right-hand side gives:
\allowdisplaybreaks
\begin{align*} \label{ub2}
    & - \left(1-C_1C_2 \sigma^{N-\tilde{N}+1}_1 \left( \frac{1-\sigma^{\tilde N }_2}{1-\sigma_2}\right) \right) \E{\P{}{}}{\ell(x_{0|k}, u_{0|k})}  \\
    &+ \sum^{N - \tilde N -1}_{j=0}\ell^0_{j} + C_2 \sigma_1 \frac{1-\sigma^{\tilde N}_2}{1-\sigma_2} \ell^0_{N-\tilde{N}} + \sum^{\tilde N}_{j=1}\ell^0_{N-\tilde{N}+j}  \\
    &\leq - \left(1-C_1C_2 \sigma^{N-\tilde{N}+1}_1 \left( \frac{1-\sigma^{\tilde N }_2}{1-\sigma_2}\right) \right) \E{\P{}{}}{\ell(x_{0|k}, u_{0|k})} \\
    & + \max\left \{1, C_2 \sigma_1 \frac{1-\sigma^{\tilde N}_2}{1-\sigma_2} \right \} \left(\sum^{N}_{j=0}\ell^0_{j} \right) \numberthis .
\end{align*}    
Comparing \eqref{ub2} with \eqref{eq:srdp_1}, we recognize 
\begin{align*}
    & \alpha = 1-C_1C_2 \sigma^{N-\tilde{N}+1}_1 \left( \frac{1-\sigma^{\tilde N}_2}{1-\sigma_2}\right), \\
    & c_2 = \max\left\{1, C_2 \sigma_1 \frac{1-\sigma^{\tilde N}_2}{1-\sigma_2} \right\} \left(\sum^{N}_{j=0}\ell^0_{j} \right). \numberthis
\end{align*}

\end{document}